\documentclass{article}

\newcommand{\sect}[1]{\section{#1}}

\newtheorem{thh}{Theorem}

\def\sep{\;\vrule\;}
\def\proof#1. {\par
                      \ifdim\lastskip<15pt
                      \removelastskip\penalty-200
                      \vskip15pt plus3pt minus3pt
                      \fi
                       {\def\a{#1}
                       \ifx\a\empty
                       {\noindent\bf Proof.}
                       \else
                       {\noindent\bf Proof of #1.}
                       \fi}\enspace}
\def\restr#1{\,\vrule\,\lower1.75ex\hbox{$#1$}}

\def\be{\begin{equation}}
\def\ee{\end{equation}}
\def\bea{\begin{eqnarray}}
\def\eea{\end{eqnarray}}
\def\bean{\begin{eqnarray*}}
\def\eean{\end{eqnarray*}}

\def\a{\alpha}
\def\b{\beta}
\def\d{\delta}

\def\g{\gamma}
\def\G{\Gamma}
\def\i{\infty}

\def\O{\Omega}

\def\s{\sigma}
\def\S{\Sigma}

\def\t{\theta}

\def\H{{\bf H}}
\def\T{{\bf T}}

\font\tenopen = cmbx10
\font\sevenopen = cmbx7
\font\fiveopen = cmbx5
\newfam\openfam
\def\open{\fam\openfam\tenopen}
\textfont\openfam = \tenopen
\scriptfont\openfam = \sevenopen
\scriptscriptfont\openfam = \fiveopen

\def\R{{\open R}}

\def\C{{\open C}}

\usepackage{graphicx}
\date{}
\title{The mean value property (corrected version)\footnote{This is a corrected version
of the paper "The mean value property, {\it Mathematical Intelligencer}, {\bf 37}(2015), 9--16.". In the
original paper an incorrect claim was stated about the connection
of harmonic functions and minimal surfaces. In the present version that incorrect claim is eliminated and
the relation in question is explained in more details in Section 8; the rest of the paper has not been changed.}}
\author{Vilmos Totik}
\begin{document}
\maketitle
\ \hskip5cm\vbox{\begin{verse}
{\footnotesize ``...what mathematics really consists of\\
is problems and solutions." (Paul Halmos)}
\end{verse}}

\sect{Some problem challenges}\label{sectintro}
I like problems and completely
agree with Paul Halmos that they are from the heart of mathematics
(\cite{Halmos}).
This one I heard when I was in high school.
\smallskip

\noindent{{\bf Problem 1.} \it Show that if numbers in between
$0$ and $1$ are written into the squares of the integer lattice on the plane in such
a way that each number is the average of the four
neighboring numbers, then
all the numbers must be the same.}
\smallskip

Little did I know at that time that this problem has
many things to do with random walks,  the fundamental
theorem of algebra, harmonic functions, the
 Dirichlet problem or with the shape of soap films.

A somewhat more difficult version is \smallskip
\smallskip

\noindent{{\bf Problem 2.} \it Prove the same if
the numbers are assumed only to be nonnegative}.\smallskip

Since 1949 every fall there is a unique mathematical contest in Hungary
named after
Mikl\'os Schweitzer, a young mathematician perished during the siege
of Budapest in 1945. It is for university students without age
groups, and about 10-12 problems from various fields
of mathematics are posted for 10 days during
which the students can use any tools and literature they
want.\footnote{The problems and solutions up to 1991
can be found in he two volumes \cite{Sch} and \cite{szekely}}
I proposed the following continuous variant of Problem 1 for
 the 1983 competition (\cite[p. 34]{szekely}).
\smallskip

\noindent{{\bf Problem 3.} \it Show that if a bounded continuous function
on the plane has the property that its average over every circle of radius
1 equals its value at the center of the circle, then it is constant.}
\smallskip

When ``boundedness" is replaced here by ``one-sided boundedness",
say positivity, the claim is
still true, but the problem gets considerably tougher.
\smallskip

\noindent{{\bf Problem 4.} \it The boundedness in Problem 3
can be replaced by positivity.}
\smallskip

We shall solve these problems and discuss
their various connections. Although the first problem follows
from the second one,  we shall first solve
Problem 1  because its solution will guide us
in the solution of the stronger statement.

It will be clear that there is nothing
special about the plane, the claims
are true in any dimension.
\begin{itemize}
\item {\it If nonnegative numbers are written into every box
of the integer lattice in $\R^d$ in such
a way that each number is the average of the $2d$
neighboring numbers,
then all the numbers are the same.}
\item {\it If a nonnegative continuous function
in $\R^d$ has the property that its average over every sphere of radius
1 equals its value at the center of the sphere, then it is constant.}
\end{itemize}

It will also be clear that similar statements are
true for other averages (like the one taken for the 9 touching
squares instead of the 4 adjacent ones).

Label a square of the integer lattice by its lower-left vertex, and
let $f(i,j)$ be the number we write into the $(i,j)$ square.
So Problem 1 asks for proving
that if $f(i,j)\in [0,1]$ and for all $i,j$
we have
\be f(i,j)=\frac{1}{4}\Bigl(f(i-1,j)+f(i+1,j)+f(i,j-1)+f(i,j+1)\Bigr),\label{11}\ee
then all $f(i,j)$ are the same.
Note that some kind of limitations like boundedness or one-sided
boundedness is needed, for, in general,  functions
with the property (\ref{11}) need not be constant, consider e.g. $f(i,j)=i$.
(\ref{11})
is called the discrete mean  value property for $f$.
First we discuss some of its consequences.

\sect{The maximum principle}\label{sectmax}
Assume that $f$ satisfies (\ref{11})
and it is nonnegative.
Notice that if, say, $f$
takes the value 0 at an $(i,j)$, then it must be zero everywhere.
Indeed, then (\ref{11}) gives that $f(i-1,j)$, $f(i+1,j)$,
$f(i,j-1)$ and $f(i,j+1)$ all must be also 0, i.e. all the
neighboring values must be zero. Repeating this we can get that
 all values of $f$ must be 0. The same argument works
 if $f$ takes its largest value at some point, so we have
 \begin{thh} \label{maxminpr}{\bf (Minimum/maximum principle)}
If a  function with the discrete mean  value property
on the integer lattice
attains somewhere its smallest/largest value, then it must be
constant.
\end{thh}
In particular, this implies a solution to Problem 1 if
we assume that $f$ has a limit at infinity (i.e. $f(i,j)\to\a$
for some $\a$ as $i^2+j^2\to\i$).
Unfortunately, in Problem 1 we do not know in advance that the
function has a limit at infinity, so this is not a solution.

\begin{figure}[h!]
\begin{center}
\includegraphics[scale=.66]{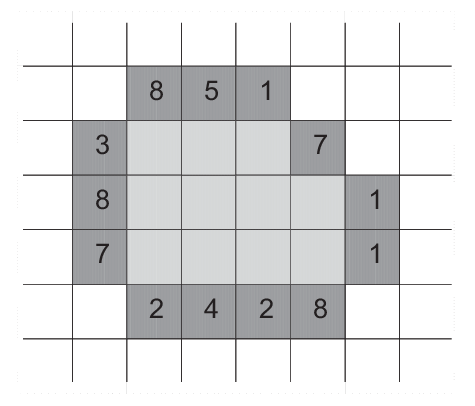}
\end{center}
\caption{\label{pict3}}
\end{figure}

Call a subset $G$ of the squares of the
integer lattice a region if every square in $G$ can be reached from every
other  square of $G$ by moving always inside $G$ to neighboring cells.
The boundary $\partial G$
of $G$ is the set of squares that are not in $G$ but which
are neighboring to $G$. See Figure \ref{pict3}
for a typical bounded region, where the boundary consists of the
darker shaded
squares.

Suppose each boundary
square contains a number like in Figure \ref{pict3}. Consider the
number filling problem:
\begin{itemize}
\item Can
the squares of $G$ be filled in with numbers so that
the discrete mean  value property is true for all squares in $G$?
\item In how many ways can such a filling be done?
\end{itemize}

This problem is called the discrete Dirichlet problem,
we shall see its connection with the classical Dirichlet
problem later.

The unicity of the solution is easy to get. Indeed,
 it is clear that the maximum/minimum principle
holds (with the proof given above) also on finite regions:
\begin{thh} \label{maxminpr1}{\bf (Maximum principle)} Let $f$ be a
function with the discrete mean  value property on a finite region,
and let $M$ be its largest value on $G\cup \partial G$.
If $f$ attains $M$ somewhere in $G$, then $f$ is the constant function.
\end{thh}
This gives that if a function with the mean value property
is zero on the boundary,
then it must be zero everywhere, and from here the unicity of the solution
to the discrete
Dirichlet problem  follows (just take the difference of two
possible solutions).

Perhaps the most natural approach to the
existence part of the
number filling problem is to consider the
numbers to be filled in as unknowns, to  write up a system of equations
for them which describes the discrete mean  value property
and the boundary properties of $f$,
and to solve that system. It can be readily shown that this
linear system of equations is always solvable.
 But there is a
better way to show existence that also works on unbounded regions.

\sect{Random walks}\label{sectrandom}
Consider a  random walk on the squares of the integer
lattice, which means that  if at
a moment we are in the lattice square
$(i,j)$, then we can move to any one of the neighboring squares
$(i-1,j)$, $(i+1,j)$, $(i,j-1)$ or $(i,j+1)$.
Which one we choose depends on some random
event, like throw two fair coins, and if the
result is ``Head-Head" then move to $(i-1,j)$,
if it is ``Head-Tail" then move to $(i+1,j)$, etc.
\begin{figure}[h!]
\begin{center}
\includegraphics[scale=.66]{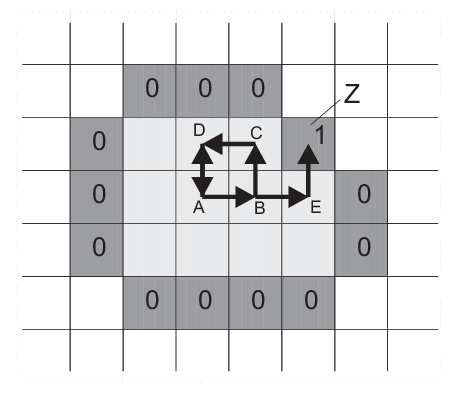}
\end{center}
\caption{A sample random walk starting at $P=A$ and terminating at
 $Z$: ABCDADABEZ \label{pict4}}
\end{figure}

We would like to find the unique
value $f(P)$ of the square-filling problem at a point $P$ of
the domain $G$. Start a random walk from $P$ which stops when
it hits the boundary of $G$. Where it stops
there is a prescribed number of the boundary, and since
it is a random event which boundary point the walk  hits first,
that boundary number is also random. Now $f(P)$ is the
expected value of that boundary number. Indeed,
from $P$ the walk  moves to either of the four neighboring squares
$P_-$, $P_+$, $P^-$ and $P^+$ with probability $1/4$--$1/4$,
and then it continues as if it was started from there.
So the just introduced expected value for $P$ will be the average of the expected
values for $P_-$, $P_+$, $P^-$ and $P^+$. Hence the
mean value property is satisfied.

Unfortunately, it is not easy to calculate the hitting probabilities
and the aforementioned expected value, but the connection with the discrete mean  value property is
notable. Furthermore, the order can be reversed, and
this connection can be used to calculate certain probabilities.
Consider the following question.
\smallskip

\noindent{\bf Problem 5}. {\it Two players, say \H\   and \T, where
\T\  is a dealer, repeatedly
 place 1--1 dollar on the table,
flip a coin, and if it is Head, then \H\  gets both
notes, while if it is Tail, then \T\  gets both of them.
Suppose \H\  starts with 30\$, and wants to know her chance
of having 100\$ at some stage, when she quits the game.}

Direct calculation of the probability of success for \H \ is non-trivial
and rather tedious.
However, using the connection between random walks and
functions with the discrete mean  value property we can easily show
that the answer is $3/10$. To this end,
let $f(i)$ be the probability of success for {\bf H} (i.e. reaching
100\$) when she
starts with $i$ dollars. After the first play {\bf H} will have
either $i-1$ or $i+1$ dollars with probability 1/2--1/2
(therefore, the fortune of  {\bf H} makes a random walk on
the integer lattice), and from
there the play goes on as if {\bf H} started with $i-1$ resp. $i+1$ dollars.
Therefore, $f$ satisfies the mean value property
\be f(i)=\frac12(f(i-1)+f(i+1)),\quad i=1,\ldots,99, \quad f(0)=0,\ f(100)=1.\label{1dim}\ee
What we have shown in dimension 2
 remains true without any change in other dimensions, in particular,
the discrete Dirichlet problem (\ref{1dim}) has one and only
one solution. But $f(i)=i/100$ is clearly a solution to (\ref{1dim}), so
 $f(30)=30/100$ as was stated
above. In general,  if {\bf H}  starts with $k$ dollars
and her goal is to reach $K$ dollars, then her chance of success is
$k/K$.
\medskip

For more on discrete random walks and their connection
with electrical circuits see the
wonderful monograph \cite{Doyle}.
We shall return to random walks in Section \ref{sectdiscretization}.

\sect{An iteration process}\label{sectiteration}
The two solutions to the discrete dirichlet problems
discussed so far (solving linear systems or
using random walks) are not too practical. Now
we discuss a fast and simple method for approximating
the solution.

Let $G$ be a bounded region and let
$f_0$ be a given function on the boundary $\partial G$.
We need to find a function $f$ on $G\cup \partial G$ which agrees
with $f_0$ on the boundary and satisfies
the discrete mean  value property in $G$.
Define for any $g$ given on $G\cup \partial G$
the function $Tg$ the following way: for $(i,j)$ in $G$
let
\be Tg(i,j)=\frac{1}{4}\Bigl(g(i-1,j)+g(i+1,j)+g(i,j-1)+g(i,j+1)\Bigr),\label{111}\ee
while on the boundary set $Tg(Q)=g(Q)$. Note that we are looking
for an $f$ that satisfies $Tf=f$,
so we are looking for a fixed point of the ``operation" $T$.
Fixed points are often found by iteration: let $g_0$ be arbitrary
and form $Tg_0,T^2g_0,\ldots$. If this happens to converge, then
the limit $f$ is a fixed point.

To start the iteration, let $g_0$ be the function
which agrees with $f_0$ on the boundary and which
is 0  on $G$.
Form $T^kg_0$, $k=1,2,\ldots$. It is a simple exercise
to show
that the iterates $T^kg_0$ converge
(necessarily to the solution of the discrete Dirichlet problem),
 and the speed of convergence is geometrically fast.

\sect{Solution to Problem 1}\label{sectsol1} Let $f$ be a function
on the lattice squares of the plain such that it has the
discrete mean  value property and its values lie in $[0,1]$.
We know from the maximum principle that if $f$ assumes somewhere
an extremal (largest or smallest) value, then $f$ is constant.

The solution uses a similar idea, by considering the set
${\cal F}$ of all such functions and considering
\be \a=\sup_{f\in{\cal F}} (f(1,0)-f(0,0)).\label{alpha}\ee
Since the translation of any $f\in {\cal F}$
by any vector $(i,j)$ is again in ${\cal F}$, and
so is the rotation of $f$ by 90 degrees, it immediately
follows that $\a$ is actually the supremum of the differences
of all possible values $f(P)-f(P')$ for neighboring squares
$P,P'$ and for $f\in {\cal F}$. Thus, Problem 1
amounts the same as showing $\a=0$ (necessarily
$\a\ge 0$).

First of all, note that there is an $f\in{\cal F}$ for which
\be \a=f(1,0)-f(0,0).\label{alpha0}\ee
Indeed, by the definition of $\a$,
for every $n$ there is an $f_n\in{\cal F}$ for which
\[ f_n(1,0)-f_n(0,0)>\a-\frac1n.\]
By selecting repeatedly subsequences we get a subsequence $f_{n_k}$
for which the sequences $\{f_{n_k}(i,j)\}_{k=1}^\i$ converge for all $(i,j)$.
 Now if
 \[f(i,j)=\lim_{k\to\i}f_{n_k}(i,j),\]
then clearly $f\in {\cal F}$ and (\ref{alpha0}) holds.

The function $g(i,j)=f(i+1,j)-f(i,j)$ has again the discrete
mean value property, and, according to what was said before,
we have $g(i,j)\le \a$ and $g(0,0)=\a$. Thus, we get
from the maximum principle that $g$ is constant, and the
constant then must be $\a$. In particular, $f(i+1,0)-f(i,0)=g(i,0)=\a$
for all $i$. Adding these for $i=0,1,,\ldots,m-1$
we obtain $f(m,0)-f(0,0)=m\a$, which is possible for
large $m$ only if $\a=0$, since $f(m,0),f(0,0)\in [0,1]$. Hence,
$\a=0$, as was claimed.

\sect{Sketch of the solution to Problem 2}\label{sectsol} Let now $f$ be a positive function
on the lattice squares of the plain such that it has the
discrete mean  value property. Without loss of
generality assume $f(0,0)=1$, and let ${\cal G}$ be the family of all such $f$'s. Then the positivity of $f$
yields that $f(0,-1),f(0,1),f(1,0),f(-1,0)\le 4$, and repeating this argument it follows that
$ 0<f(i,j)\le 4^{|i|+|j|}$ for all $f\in {\cal G}$ and for all $i,j$.
Hence, the selection process in the preceding section can be carried out
without any change in the family ${\cal G}$.
Now  consider
\be \b:=\sup_{f\in {\cal G}}f(1,0)=\sup_{f\in {\cal G}}\frac{f(1,0)}{f(0,0)}.
\label{beta}\ee
As before, $\b$ turns out to be the supremum of the ratios $f(P')/f(P)$
for all neighboring squares $P,P'$ and for all $f\in {\cal G}$, therefore
   Problem 2 asks for showing that
$\b=1$ (clearly $\b\ge 1$).

Let $f\in {\cal G}$ be a function for which equality is assumed
in (\ref{beta}) (the existence of $f$ follows from the selection process).
Then
\be \b=\frac{f(1,0)}{f(0,0)}=
\frac{ f(2,0)+f(1,1)
+f(1,-1)+ f(0,0)}
{f(1,0)+f(0,1)
+f(0,-1)+ f(-1,0)},\label{al0}\ee
which
can only be true if each one of the upper terms equals $\b$ times the term
below it (since each term in the numerator is at most $\b$ times the term right below it).
In particular, $f(2,0)/f(1,0)=\b$ and $f(0,0)/f(-1,0)=\b$. Repeat the previous argument
to conclude that
$f(k,0)/f(k-1,0)=\b$ for all $k=0,\pm1,\pm2,\ldots$. Thus,
there is some constant $\g_0>0$ such that $f(k,0)=\g_0 \b^k$
for all $k$.
Since
$f(1,1)/f(0,1)=\b$ and $f(1,-1)/f(0,-1)=\b$
are also true, it follows as before that
$f(k,\pm1)=\g_{\pm1}\b^k$ for all $k$ with some $\g_{\pm 1}>0$.
Repeating again this
argument we finally conclude that there are positive numbers
$\g_j$ such that $f(i,j)=\g_j\b^i$ is true for all $i,j$.

Apply now the discrete mean  value property:
\be \g_j=f(0,j)
=\frac{1}{4}\Bigl(\g_j\frac1\b+\g_j\b +\g_{j-1}+\g_{j+1}\Bigr).\label{nm}\ee
Since $\b+1/\b\ge 2$, this implies $2\g_j\ge \g_{j-1}+\g_{j+1}$,
i.e. the concavity of the sequence $\{\g_j\}_{-\i}^\i$. But a positive sequence
on the integers can be concave only if it is constant.
Thus, all the $\g_j$'s are the same, and then (\ref{nm})
cannot be true if $\b>1$, hence $\b=1$ as claimed.

\sect{The continuous mean value property and harmonic functions}\label{sectcont}
Assume that $G\subset \R^2$ is a domain
(a connected open set) on the plane, and $f:G\to \R$ is
a continuous real-valued function defined on $G$. We say that $f$ has
the mean value property in $G$ if for every circle
$C$
which lies in $G$ together with its interior we have
\be f(P)=\frac1{|C|}\int _C f,\label{cont1}\ee
where $P$ is the center of $C$ and $|C|$ denotes the length of $C$.
(\ref{cont1}) means that the average of $f$ over the circle $C$
coincides with the function value at the center of $C$.

Functions with this mean value property are called harmonic,
and they play a fundamental role in mathematical analysis.
For example, if $f$ is the real part of a complex differentiable
(so called analytic) function, then $f$ is harmonic.  The converse is also
true in simply connected domains (domains without holes): if $f$ has the mean
value property (harmonic), then it is the real part
of a differentiable complex function. So there is an
abundance of harmonic functions, e.g. the real part of
any polynomial is harmonic, say
\[f(x,y)=\Re z^n=x^n-{n\choose 2} x^{n-2}y^2+{n\choose 4} x^{n-4}y^4-\cdots\]
are all harmonic.

Although we shall not use it, we mention that the standard (but equivalent)
definition of harmonicity is $f_{xx}+f_{yy}=0$, where
$f_{xx}$ and $f_{yy}$ denote the second partial derivatives of
$f$ with respect to $x$ and $y$. We shall stay with our geometric definition.

Simple consequence of the mean value property is the maximum principle:
\begin{thh} \label{maxminpr10}{\bf (Maximum principle)}
If a harmonic  function on a domain $G$
attains somewhere its largest value, then it must be
constant.
\end{thh}
The reader can easily modify the argument given for
Theorem \ref{maxminpr} to verify this version.

A basic fact concerning harmonic functions is
that a bounded harmonic function on the
whole plane must be constant (Liouville's theorem). I heard the following
proof from Paul Halmos.
Suppose $f$ satisfies (\ref{cont1}) on the whole
plane and it is bounded. Let $D_R(P)$ be the disk of radius $R$ about some
point $P$. Since
the integral over $D_R(P)$ can be obtained by first integrating
on circles $C_r$ of radius $r$ about $P$ and then integrating
 these integrals with respect to $r$
(from 0 to $R$), it easily follows that $f$ also has
the area-mean value property:
\be f(P)=\frac1{{R^2\pi}}\int _{D_R(P)} f.\label{cont2}\ee
If $Q$ is another points, then the same
formula holds for $f(Q)$ with $D_R(P)$ replaced by $D_R(Q)$.
Now for very large $R$ the disks $D_R(P)$ and $D_R(Q)$
are ``almost the same" in the sense that outside their
common part there are only two small regions in them
the area of which is negligible compared to the area
of the disks. So the averages of $f$ over $D_R(P)$ and
$D_R(Q)$ are practically the same (by the boundedness of $f$),
and for $R\to\i$ we get that in the limit the
averages, and hence also the function values
at $P$ and $Q$ are the same.

From here the fundamental theorem of algebra (``every polynomial
 has a zero on the complex plane") is a standard consequence
(if the polynomial $\cal P$ did not vanish anywhere, then the real and imaginary
parts of $1/\cal P$ would be bounded harmonic functions on the plane, hence they
would be constant, which is not the case).

Note that Problem 3 claims more than Liouville's theorem,
since in it the mean
value property is requested only for circles with
a fixed radius. In general, if we know the
mean value property (\ref{cont1}) for all circles
of a fixed radius $C=C_{r_0}$, then it does not
follow that $f$ is harmonic. However, by a result of Jean Delsarte
if (\ref{cont1}) is true for all circles
of radii equal to some $r_0$ or $r_1$ and
$r_0/r_1$ does not lie in a finite
exceptional set (consisting of the ratios of solutions of an
equation involving a Bessel function), then $f$ must
be harmonic. In $\R^3$ this exceptional set is empty,
and it is conjectured that it is empty in all
$\R^d$, $d> 2$. See the most interesting
paper \cite{Zalcman}, as well as the extended literature
on Pompeiu's problem in \cite{Zalcman1}--\cite{Zalcman2}.

\sect{The Dirichlet problem and soap films}\label{sectdirichlet}
The Dirichlet problem in the continuous case is the following:
suppose $\O $ is a (bounded) domain with boundary $\partial \O $,
and there is a continuous function $g_0$ given on the
boundary $\partial \O $. Can this $g_0$ be (continuously) extended
to $\O$ to a harmonic function, i.e. can it be continued inside
$\O$ so that it has the mean value property there?
Under some minimal regularity assumptions on the boundary
this problem has always a unique solution $g$.

\begin{figure}[h!]
\begin{center}
\includegraphics[scale=.5]{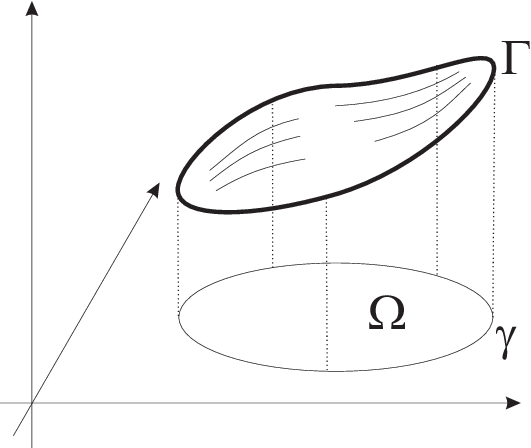}
\end{center}
\caption{The plane curve $\g$ and its lift-up $\G$\label{pict6}}
\end{figure}

The unicity
follows from the maximum principle in the same fashion
that was done in the discrete case. The existence requires additional assumptions,
for example if $\O $ is the punctured disk $\{z\sep 0<|z|<1\}$ and
we set $g_0(z)=0$ for $|z|=1$ while $g_0(0)=1$, then there is
no harmonic function in $\O $ which continuously extends $g_0$.
But this example is pathological
($0$ is an isolated point on the boundary), and it can be
shown that in most cases the Dirichlet problem can be solved.
Later we sketch how.

Suppose that the boundary of the domain $\O $
is a simple closed curve $\g$
with parametrization $\g(t)=(\g_1(t),\g_2(t))$.
Consider the given function $g_0$ on $\partial \O =\g$,
and with its help lift  $\g$ up into 3 dimensions:
$\G(t):=(\g_1(t),\g_2(t),g_0(\g(t)))$ is a 3 dimensional curve above
$\g$. Now the points
$(x,y,g(x,y))$, $(x,y)\in \O$, describe a surface that is
attached to that space curve $\G$.

There are other surfaces attached to $\G$, among the most studied ones
is the so called minimal surface with boundary $\G$.
$\G$ can be thought of as
 a wire, and stretch an elastic rubber sheet
(or a soap film) over $\G$ (see Figure \ref{pict6}). When in rest,
the rubber sheet (soap film) gives a surface ${\cal H}_\G$ over the domain $\O $,
which is the graph of a function $f$. In general, $f$ is not the
same as $g$ (the solution of the Dirichlet problem), i.e.
$f$ is not harmonic (it can be characterized by a partial
diffential equation). Nonetheless, ${\cal H}_\G$ (a so-called minimal surface -- minimizing
the area) can also be characterized via harmonic functions.
To do that we need to introduced isothermal parametrizations.

Suppose we are given three continuous functions
$x_1(u,v),x_2(u,v),x_3(u,v)$ defined in a domain $D$ on the plane.
Under some simple regularity conditions the points $(x_1(u,v),x_2(u,v),x_3(u,v))$, $(u,v)\in D$,
describe a surface ${\cal S}$ in $\R^3$. In this case we say that $(x_1(u,v),x_2(u,v),x_3(u,v))$, $(u,v)\in D$,
gives ${\cal S}$
in parametric form. If $(u(t),v(t))$, $t\in [0,1]$, is
a curve $\s$ in $D$, then
\[(x_1(u(t),v(t)),x_2(u(t),v(t)),x_3(u(t),v(t))), \qquad t\in [0,1],\]
is a space curve $\S$ that lies in ${\cal S}$. There are different ways to
parametrize a surface ${\cal S}$, and a parametrization is called isothermal if
it preserves the angle between curves, meaning that  if $\s_1$ and
$\s_2$ are (smooth) curves in $D$ that intersect each other, then the
angle between $\s_1$ and $\s_2$ at their intersection point\footnote{The angle between two curves
intersecting each other in a point $P$ is the angle that the tangents to the curves at $P$ form.}
is the same as the angle between the corresponding space curves $\S_1$ and $\S_2$ at their
corresponding intersection point. Under some regularity conditions such an isothermal parametrization
exists (at least locally, c.f. \cite[p. 31]{Osserman}.

Now it turns out (see e.g. \cite[Lemma 4.2]{Osserman}) that a surface ${\cal S}$ 
attached to ${\bf \G}$ with isothermal parametrization is the  minimal surface
${\cal H}_\G$ if and only if all the parameter functions $x_1(u,v)$, $x_2(u,v)$, $x_3(u,v)$
in its isothermal parametrization are harmonic funtions.

\sect{Discretization and Brownian motions}\label{sectdiscretization}
Is there a connection between the discrete and continuous Dirichlet
problems discussed in Sections \ref{sectmax} and \ref{sectdirichlet}?
There is indeed, and it is of practical importance.

\begin{figure}[h!]
\begin{center}
\includegraphics[scale=.66]{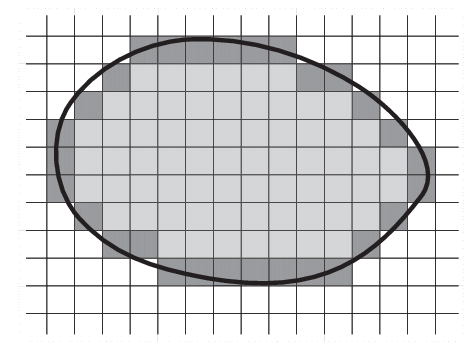}
\end{center}
\caption{The domain $\O$ enclosed by the closed curve and
the region $G$ of squares lying inside $\O$
(with darker shaded boundary squares)\label{pict5}}
\end{figure}

Let $\O$ be a domain as in the preceding section,
and let $g_0$ be a continuous function on the boundary of $\O$.
Consider a square lattice on the plane with small mesh size, say
consisting of $\tau\times\tau$ size squares with small $\tau$.
Form a region $G_\tau$ (see Figure \ref{pict5})
on this lattice by considering those
squares in the lattice which lie in $\O$ (there may
be a slight technical trouble that the union of these squares
may not be connected, in that case let $G_\tau$ be the union
of all squares that can be reached from a square containing
a fixed point of $\O$). We are going to consider
the discrete Dirichlet problem on $G_\tau$, the solution of which
will be close to the solution of the original continuous Dirichlet problem.
To this end define a boundary function on the boundary
squares $\partial G_\tau$ in our lattice: if $P$ is a boundary square,
then $P$ must intersect the boundary $\partial \O$, and if
$z\in P\cap \partial \O$ is any point, then set $f_{0,\tau}(P)=g_0(z)$.
Now solve this discrete Dirichlet problem on $G_\tau$ (with boundary
numbers given by $f_{0,\tau}$) with the iteration
technique of Section \ref{sectiteration}. Note that the iteration
in Section \ref{sectiteration} is computationally very simple and quite
fast, since all one needs to do is to calculate averages of 4--4 numbers.
Besides that,  the convergence of the iterants
to the solution is geometrically fast.
Let $f_\tau$ be the solution, and we can imagine that $f_\tau$ gives
us a function $F_\tau$ on the union of the squares belonging to $G_\tau$:
on every square $P$ the value of this $F_\tau$ is identically equal
to  the number $f_\tau(P)$. Now if $\tau$ is small, then this function
$F_\tau$ will be close on $G_\tau$ to the solution $g$ of the continuous Dirichlet
problem we are looking for.

There is yet another connection between the discrete
and the continuous Dirichlet problems.
We have seen in Section \ref{sectrandom} that the discrete Dirichlet problem
can be solved via random walks on the squares
of the integer lattice. Now consider the just-introduced
square  lattice with small mesh size, and make a random
walk on that lattice. If the mesh size is getting smaller then
the lattice is getting
denser (alternatively look at the square lattice from a far distance).
To compensate for having more and more squares,  speed up the random walk.
If this speeding-up is done properly, then in the limit we
get a random motion  on the plane, the Brownian motion.
In a Brownian motion a particle moves in such a way that it continuously and randomly
changes its direction.
\begin{figure}[h!]
\begin{center}
\includegraphics[scale=.66]{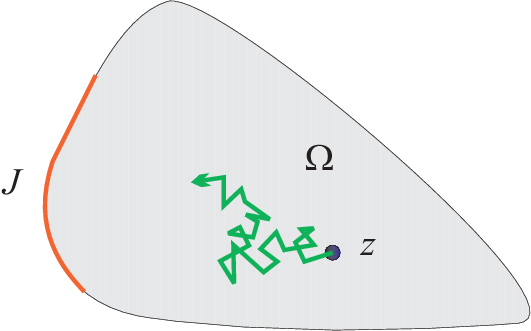}
\end{center}
\caption{A Brownian motion\label{pict7}}
\end{figure}

Let $\O$ have smooth boundary,
and let $J$ be an arc on that boundary, see Figure \ref{pict7}. Start a
Brownian motion at a point $z\in \O $, and stop it when it
hits the boundary of $\O $, and  let
$f_J(z)$ be the probability that it hits the boundary
 in a point of $J$. This $f_J(z)$ (which is called
 the harmonic measure of $z$ with respect to $\O$ and $J$)
 has the mean value property.
Indeed, consider a circle $C$ about the point $z$ that lies inside
$\O $ together with its interior. During the motion of the particle
there is a fist time when the particle hits $C$ at a point $Z\in C$. Then it
continues as if it started in $Z$, and then
the probability that it hits the boundary $\partial \O $
in a point of $J$ is $f_J(Z)$. Because of the circular symmetry
of $C$, all $Z\in C$ play equal roles, and we can conclude
(at least heuristically), that the hitting probability
$f_J(z)$ is the average of the hitting probabilities
$f_J(Z)$, $Z\in C$, which is precisely the mean value property
for $f_J$.
It is also clear that if $z\in \O$ is close to a point
$Q$ on the boundary of
$\O$, then  it is likely that the Brownian motion starting in
$z$ will hit the boundary $\partial \O $ close to $Q$. Therefore, if
$Q\in J$ (except when  $Q$ is one of the endpoints of $J$) the probability
$f_J(z)$ gets higher and higher, eventually converging to
1 as $z\to Q$, while in the case when $Q\not\in J$,
the probability $f_J(z)$ gets smaller and smaller,
eventually converging to 0.

What we have shown is that $f_J$ is a harmonic function in $\O $
which extends continuously to the boundary  to 1 on the inner part of the
arc $J$ and it extends continuously to 0 on the outer
part of $J$ (therefore, at the endpoints $f_J$ cannot  have
a continuous extension). In other words, not worrying about
continuity at the endpoints of the arc $J$, we have
solved the Dirichlet problem for the characteristic function
\[\chi_J(z)=\left\{\begin{array}{ll}
  1 & \mbox{if $z\in J$} \\
  0 & \mbox{if $z\not\in J.$}
\end{array} \right.\]

Now if $f_0$ is a continuous function
on the boundary $\partial \O $,
then to $f_0$ there is arbitrarily close a function of the
form $h=\sum c_j \chi_{J_j}$ with a finite sum, and then $f_h:=\sum c_jf_{J_j}(z)$
is a harmonic function in $\O $ which is close to $f_0$ on the boundary.
Using the maximum principle it follows that, as $h\to f_0$,
the functions $f_h$ converge uniformly on $\O \cup \partial \O $
to a function $f$ which has the mean value property in
$\O $ (since all $f_h$ had it) and which agrees with $f_0$
on $\partial \O $. Therefore, this $f$ solves the Dirichlet problem.
\medskip

We refer the interested reader to the book \cite{Ransford} for
further reading concerning the mean value property and the Dirichlet problem.
Robert Brown in ``Brownian motion" was a Scottish
botanist who, in 1827, observed in a microscope that pollen particles in
suspension make an irregular, zigzag motion.
The rigorous mathematical foundation of Brownian motion
was made by Norbert Wiener in 1923 (\cite{Wiener}).
The connection to the Dirichlet problem was first observed by
Shizuo Kakutani \cite{Kakutani}. This had
a huge impact on further developments;
there are many works that discuss the relation between random walks
and problems (like the Dirichlet problem) in potential theory,
see e.g. \cite{Port} or the very extensive \cite{Doob}.

\sect{Solution to Problem 4}
In this proof we shall be brief, since
some of the arguments have already
been met before.

First of all, seemingly nothing prevents an $f$ as in Problem 4
behave wildly, and first we ``tame" these functions.
Let $\cal F$ be the collection of all positive functions
on the plane with the mean value property (\ref{cont1}),
and for some $\d>0$ let ${\cal F}_\d$ be the collection
of all the functions
\be f_\d(z)=\frac{1}{\d^2 \pi}\int_{D_\d(z)} f(u)du\label{av1}\ee
 for $f\in {\cal F}$, where $D_\d(z)$ denotes the disk
 of radius $\d$ about the point $z$.
If we can show that
\[\b_\d:=\sup_{g\in {\cal F}_\d,\ z\in \C, \ |\t|=1} g(z+\t)/g(z)\]
is 1, then we are done. Indeed, then $g(z+\t)\le g(z)$ holds
 for all $g\in {\cal F}_\d$,
$z\in \C$ and any $\t$ with $|\t|=1$, which actually implies
 $g(z+\t)= g(z)$ (just apply the inequality to $z+\t$ and to
 $-\t$). Hence, since any two points on the plane can be connected
 by a polygonal line consisting of segments of length 1,
  every $g=f_\d\in {\cal F}_\d$ is constant, and then
letting $\d\to 0$ we get that every $f$ in ${\cal F}$ is constant,
as Problem 4 claims.

First we show that $\b_\d$ is finite. From the mean value property
(\ref{cont1}) we have for $f\in {\cal F}$
\bea f(z)=\frac{1}{2\pi}\int_0^{2\pi}f(z+e^{it})dt&=&
\frac{1}{2\pi}\int_0^{2\pi}\frac{1}{2\pi}\int_0^{2\pi}f(z+e^{it}+e^{iu})dudt
\nonumber \\
&=&\int_{D_2(0)}f(z+w)A(w)d|w|,\label{jj}\eea
where $D_2(0)$ is the disk of radius 2 about the origin,
$d|w|$ denotes area-integral and $A(w)$ is a function on
$D_2(0)$ that is continuous and positive for $0<|w|<1$.
Thus, if $S$ is the ring $1/2\le |w|\le 3/2$ and $a>0$ is the minimum
of $A(w)$ on that ring, then
\[f(z)\ge a \int _{z+S} f.\]
Since for $|\t|=1$ and for $|z'-z|\le \d\le 1/4$ the ring $z'+S$ contains the disk
$D_\d(z+\t)$, it follows that
\[f(z')\ge a \int _{D_\d(z+\t)} f,\]
and upon taking the average for $z'\in D_\d(z)$, the inequality $f_\d(z)\ge (a\d^2\pi) f_\d(z+\t)$ follows.
Hence, $\b_\d$ is finite.

From the finiteness of $\b_\d$ it follows that if $R>0$ is given
and $|u|\le R$, then
$a_R\le g(z+u)/g(z)\le A_R$
for all $g\in {\cal F}_\d$ and all $z\in \C$,
where the constants
$a_R,A_R>0$ depend only on $R$.
Let  ${\cal F}_{\d\d}$ be the collection of all $g_\d$
with $g \in {\cal F}_{\d}$. Then  ${\cal F}_{\d\d}\subseteq {\cal F}_\d$,
and $\b_\d=1$ follows from $\b=1$ (to be proven
in a moment), where
\[\b:=\sup_{h\in {\cal F}_{\d\d},\ z\in \C, \ |\t|=1} h(z+\t)/h(z).\]
Since ${\cal F}_{\d\d}$ is translation- and rotation-invariant,
it is clear that
\be 1\le \b=\sup_{h\in {\cal F}_{\d\d}, h(0)=1} h(1).\label{ggh}\ee
But the collection $h\in {\cal F}_{\d\d}$
with $h(0)=1$ consists of functions that are uniformly bounded and
uniformly equicontinuous on all disks $D_R(0)$, $R>0$, hence
from every sequence of such functions one can select a
subsequence that converges uniformly on all the disks
$D_R(0)$, $R>0$,.
Therefore, the supremum in (\ref{ggh}) is attained, and
there is an extremal function $h\in {\cal F}_{\d\d}$ with
$h(1)=\b h(0)$. Suppose that $h(z+1)=\b h(z)$ holds for some
$z$ (we have just seen that $z=0$ is such a value).
From (\ref{jj}) it follows then that
\[\int_{D_2(0)}h(z+1+w)A(w)d|w|=\b\int_{D_2(0)}h(z+w)A(w)d|w|,\]
and since here, by the definition of $\b$,
$h(z+1+w)\le \b h(z+w)$ for all $w$, we can conclude
that $h(z+1+w)= \b h(z+w)$ must be true for all $|w|\le 2$.
Thus, $h(1)=\b h(0)$ implies $h(w+1)=\b h(w)$ for all $|w|\le 2$,
and repeated application of this step gives
that $h(z+1)=\b h(z)$ holds for all $z$.

Now let ${\cal F'}$ be the collection of all $h\in {\cal F}_{\d\d}$ that satisfies
the just established functional equation $h(z+1)=\b h(z)$, and let
\[ \g:=\sup_{f\in {\cal F}',\ z\in \C} h(z+i)/h(z).\]
Since ${\cal F}'$ is closed for translation and the operation $z\to \overline z$
(complex conjugation) taken in the argument,
it follows that $\g\ge 1$,
and the reasoning we just gave for $\b$ yields that
there is an extremal function $h'\in {\cal F}'$ such that
$h'(i)=\g h'(0)$, and for this extremal function we have
the functional equation $h'(z+i)=\g h'(z)$ for all $z\in \C$. Thus, $h'$
satisfies both equations $h'(z+1)=\b h'(z)$, $h'(z+i)=\g h'(z)$, from which
it follows that if $m$ is the minimum
of $h'$ on the unit square, then $h'(z)\ge m \b^i \g^j$
at the integer lattice cell with lower left corner at $(i,j)$.

Suppose now to the contrary that $\b>1$. Then the preceding estimate
shows that $h'(z)\to \i$ as the real part of $z$ tends to infinity and the imaginary
part stays nonnegative (then $i\to\i$, $j\ge 0$).
Thus, if $h''(z)=h'(z)+ h'(\overline z)$,
then $h''(z)\to\i$ as the real part of $z$ tends to $\i$,
and then
\[h'''(z)=h''(z)+h''(zi)+h''(zi^2)+h''(zi^3)\]
is a function with the mean value property which
tends to infinity as $z\to\i$ (note that $h'$ is positive,
so $h''$, $h'''$ are larger than any of the terms on the right
of their definitions). But this contradicts the maximum/minimum principle,
and that contradiction proves that, indeed, $\b=1$.

\sect{The Krein-Milman theorem}
Although all the proofs we gave were elementary,
one should be aware of a general principle about extremal
points that lies behind
these problems. Recall that
in a linear space a point $P\in K$ is called an extremal point
of a convex set $K$ if $P$ does not lie inside any segment joining
two points of $K$.

A linear topological space is called locally convex
if the origin has a neighborhood basis consisting of
convex sets.
For example, $L^p$-spaces are
locally convex precisely for $p\ge 1$. Now a theorem
of Mark Krein and David Milman says that if $K$ is a compact
convex set in a locally compact topological space, then
$K$ is the closure of the convex hull of its
extremal points.

A point $P$ is an extremal point for a convex set $K$ precisely if it has
the property that if $P$ lies in the convex hull
of a set $S\subset K$, then $P$ must be one of the points of $S$.
Now functions with a mean value property similar to
those we considered in this article form a convex set $K$ (in the
locally convex topological space of continuous or discrete functions), and the mean value property
itself means that each such function lies in the
convex hull of some of its translates. Therefore, such a function
can be an extremal point for $K$ only if it agrees with all those translates,
which means that it is constant. Now if the extremal points in $K$
are constants, then so are all functions in $K$ provided
we can apply the Krein-Milman theorem. Hence, the
crux of the matter is to prove that the additional
boundedness or one-sided boundedness hypotheses set forth in our
problems imply
that $K$ is compact; then the Krein-Milman theorem
finishes the job.
In our proofs we faced the same
problem: we needed the existence of extremal functions in
(\ref{alpha}), (\ref{beta}), (\ref{ggh}), for which we needed to prove
some kind of compactness.

In conclusion we mentioned that the problems that have been discussed
in this paper are special cases of the Choquet-Deny convolution
equation first discussed by Gustave Choquet and Jaques Deny in 1960,
which has applications in probability theory and far reaching generalizations
in various groups/spaces.
See \cite{Deny1}, \cite{Deny2} and
the extended list of references  in \cite{chu}.

\bigskip

{\bf Acknowledgement.} Part of this work was presented
in 2011 at the Fazekas Mih\'aly Gimn\'azium, Budapest, Hungary
as a lecture for high school students.
The author is thankful to Andr\'as Hrask\'o and J\'anos Pataki
for their comments regarding the presentation.
The preparation of this paper was partially supported by NSF DMS-1265375.
The author also thanks the anonymous referee for her/his helpful
suggestions.

Special thanks go to Jos\'e Gonz\'ales Llorente for pointing out the
error in the original paper about the connection of harmonic functions and
minimal surfaces.

Bolyai Institute

University of Szeged

Szeged

Aradi v. tere 1, 6720, Hungary

\smallskip

{\it totik@mail.usf.edu}


\begin{thebibliography}{9}
\bibitem{Deny1} G. Choquet and J. Deny, Sur l\'{}\'equation de convolution $\mu=\mu*\sigma$,
{\it  C.R. Acad. Sc. Paris}, {\bf 250}(1960), 779--801.

\bibitem{chu} C. Chu and T. Hilberdink, The convolution equation
of Choquet and Deny on nilpotent groups, {\it Integr. Equat. Oper. Th.},
{\bf 26}(1996), 1--13.

\bibitem{Sch} {\it Contests in Higher Mathematics, 1949-1961}, Akad\'emiai Kiad\'o, Budapest,
1968.

\bibitem{Deny2} J. Deny, Sur l\'{}\'equation de convolution $\mu=\mu*\sigma$,
{\it S\'emin. Th\'eor. Potentiel de M. Brelot},
Paris 1960.

\bibitem{Doob} J. L. Doob, {\it Classical potential theory and its probabilistic counterpart},
Grundlehren der Mathematischen Wissenschaften, {\bf 262},  Springer-Verlag, New York, 1984.

\bibitem{Doyle} P. G. Doyle and J. L.  Snell, {\it Random walks and electric networks},
Carus Mathematical Monographs, {\bf 22}. Mathematical Association of America, Washington, DC, 1984.



\bibitem{Halmos} P. Halmos, The heart of mathematics,
{\it Amer. Math. Monthly}, {\bf 87}(1980), 519--524.

\bibitem{Kakutani} S. Kakutani, Two-dimensional Brownian motion and harmonic functions,
{\it   Proc. Imp. Acad.}, {\bf  20}(1944), 706--714.


\bibitem{Osserman} R. Osserman, {\it A survey of minimal surfaces,} Van Nostrand Reinhold Mathematical Studies,
{\bf 25}, New York, 1969.

\bibitem{Port} S. C. Port and C. J. Stone, {\it Brownian motion and classical potential theory},
Probability and Mathematical Statistics, Academic Press, New York-London, 1978.

\bibitem{Ransford} T. Ransford, {\it Potential Theory in the
Complex Plane},
Cambridge University Press, Cambridge, 1995

\bibitem{szekely} G. Sz\'ekely (editor), {\it Contests in Higher Mathematics},
Problem Books in Mathematics, Springer Verlag, New York, 1995.


\bibitem{Wiener} N. Wiener, Differential space, {\it Journal of Mathematical Physics}, {\bf  2}(1923) 131--174.



\bibitem{Zalcman} L. Zalcman, Offbeat integral geometry, {\it Amer. Math. Monthly},
{\bf 87}(1980), 161--175.

\bibitem{Zalcman1} L. Zalcman, A Bibliographic Survey of the Pompeiu Problem,
Approximation by Solutions of Partial Differential Equations (Hanstholm, 1991),
{\it NATO Adv. Sci. Inst. Ser. C Math. Phys. Sci.}, {\bf 365},
Kluwer Acad. Publ., Dordrecht, 1992, 185--194.

\bibitem{Zalcman2} L. Zalcman, Supplementary bibliography to: "A bibliographic survey of the Pompeiu problem'' [in Approximation by solutions of partial differential equations (Hanstholm, 1991), 185--194, Kluwer Acad. Publ., Dordrecht, 1992],
Radon transforms and tomography (South Hadley, MA, 2000),  {\it Contemp. Math.}, {\bf 278},
Amer. Math. Soc., Providence, RI, 2001, 69--74.


\end{thebibliography}
\end{document}